\documentclass[10pt]{article}

\begin{document}

\author{Hans-Peter A. K\"unzi and Dominic van der Zypen}
\title{Maximal (sequentially) compact topologies}

\date{
{\small Dedicated to Professor Horst Herrlich on the occasion of
his $65^{th}$ birthday}}

\newtheorem{Prp}{\qquad Proposition}
\newtheorem{Thm}{\qquad Theorem}
\newtheorem{Lem}{\qquad Lemma}
\newtheorem{Cor}{\qquad Corollary}
\newtheorem{Rem}{\qquad Remark}
\newtheorem{Ex}{\qquad Example}
\newtheorem{Prb}{\qquad Problem}
\newtheorem{Def}{\qquad Definition}
\maketitle

\begin{abstract}
\noindent We revisit the known problem whether each compact
topology is contained in a maximal compact topology and collect
some partial answers to this question. For instance we show that
each compact topology is contained in a compact topology in which
convergent sequences have  unique limits. We also answer a
question of D.E. Cameron by showing that each sequentially compact
topology is contained in a maximal sequentially compact topology.
We finally observe that each sober compact $T_1$-topology is
contained in a maximal compact topology and that each  sober
compact $T_1$-topology which is locally compact or sequential is
the infimum of a family of maximal compact topologies.

\noindent \footnotetext{\noindent
AMS (2000) Subject Classifications: 54A10, 54B15, 54D10, 54D30, 54G20 

The first author acknowledges support under the bilateral
cooperation betweens Flanders and South Africa (period 2003/4).

 }
\end{abstract}

\section{Introduction}

A topological space is called a $KC$-{\em space} (compare also
\cite{Fleissner}) provided that each compact set is closed. A
topological space is called a $US$-{\em space} provided that each
convergent sequence has a unique limit. It is known
\cite{Wilansky} that each Hausdorff space (= $T_2$-space)
is a $KC$-space, each
$KC$-space is a $US$-space and each $US$-space is a $T_1$-space
(that is, singletons are closed); and no converse implication
holds, but each first-countable $US$-space is a Hausdorff space.

A compact topology on a set $X$ is called {\em maximal compact}
provided that it is not strictly contained in a compact topology
on $X.$ It is known that a topological space is maximal compact if
and only if it is a $KC$-space that is also compact
\cite{Ramanathan}. These spaces will be called compact $KC$-spaces
in the following.

Let us note that while there are many maximal compact topologies,
minimal noncompact topologies do not exist: Any noncompact space
$X$ possesses a strictly increasing open cover $\{C_\alpha:\alpha
<\delta\}$ of $X$ where $\delta$ is a limit ordinal and $C_0$ can
be assumed to be nonempty. Clearly then $\{\emptyset,X\}\cup
\{C_\alpha:0<\alpha <\delta\}$ yields a base of a strictly coarser
noncompact topology on $X.$

 Maximal compact topologies need not be
Hausdorff topologies \cite{HingTong} (see also
\cite{Balachandran,SmytheWilkins}). A standard example of a
maximal compact topology that is not a Hausdorff topology is given
by the one-point-compactification of the set of rationals equipped
with its usual topology.

Indeed maximal compact spaces can be anti-Hausdorff (=
irreducible), as we shall next observe by citing an example due to
van Douwen (see \cite{vanDouwen}).

In order to discuss that example we first recall some pertinent
definitions. A nonempty subspace $S$ of a topological space is
called {\em irreducible} (see e.g. \cite{Hoffmann}) if each pair
of nonempty open sets of $S$ intersects. Furthermore a topological
space $X$ called a {\em Fr\'echet space} (see \cite[p.
53]{Engelking}) provided that for every $A\subseteq X$ and every
$x\in \overline{A}$ there exists a sequence of points of $A$
converging to $x.$ For the convenience of the reader we include a
proof of the following observation (compare e.g. Math. Reviews
53\#1519 of \cite{Review}).

\begin{Lem} Each Fr\'echet $US$-space $X$ is a $KC$-space.
\end{Lem}

{\em Proof.} Suppose that $x\in \overline{K}$ where $K$ is a
compact subspace of $X.$ Because $X$ is a Fr\'echet space, there
is a sequence $(k_n)_{n\in {\bf N}}$ of points of $K$ converging
to $x.$ Since $K$ is compact, that sequence has a cluster point
$c$ in $K.$ Because $X$ is a Fr\'echet space, there is a
subsequence of $(k_n)_{n\in {\bf N}}$ converging to $c$ (compare
\cite[Exercise 1.6D]{Engelking}). Hence $x=c\in K,$ because $X$ is
a $US$-space. We have shown that $K$ is closed and conclude that
$X$ is a $KC$-space.

\begin{Ex} (van Douwen \cite{vanDouwen}) There exists a countably
infinite compact Fr\'echet $US$-space that is anti-Hausdorff. By
the preceding lemma that space is a $KC$-space and hence maximal
compact. Thus there exists an infinite maximal compact space that
is irreducible.
\end{Ex}

On the other hand, by the result cited
above each first-countable maximal compact ($T_1$-)topology
satisfies the Hausdorff condition (compare \cite[Theorem
8]{Stone}).

\section{Main problem and related questions}

While it is known that each compact topology is contained in a
compact $T_1$-topology (just take the supremum of the given
topology with the cofinite topology) \cite[Theorem 10]{Stone}, the
question whether each compact topology is contained in a compact
$KC$-topology (that is, is contained in a maximal compact
topology) seems still to be open. Apparently that question was
first asked by Cameron \cite[p. 56, Question 5-1]{Cameron}, but
remained unanswered.

Of course, a simple application of Zorn's Lemma cannot help us
here, since a chain of compact topologies need not have a compact
supremum: Consider the sequence $(\tau_n)_{n\in {\bf N}}$ of
topologies $\tau_n=\{\emptyset,{\bf N}\}\cup \{[1,k]:k\in {\bf
N},k\leq n\}$ $(n\in {\bf N})$ on the set ${\bf N}$ of positive
integers.

On the other hand, for instance each infinite topological space
$X$ with a point $x$ possessing only cofinite neighborhoods is
clearly contained in a maximal compact topology: Just consider the
one-point-compactification $X_x$ of $X\setminus \{x\}$ where
$X\setminus \{x\}$ is equipped with the discrete topology and $x$
acts as the point at infinity.

The problem formulated above seems to be undecided even under
additional strong conditions. Recall that a topological space is
called {\em locally compact} provided that each of its points has
a neighborhood base consisting of compact sets.
Note that a locally compact $KC$-space is a regular Hausdorff space.

\begin{Prb} Is each locally compact (resp. second-countable) compact
topology contained
in a maximal compact topology?
\end{Prb}

The authors also do not know the answer to the following
generalization of their main problem.

\begin{Prb} Is each compact topology the continuous image of a maximal
compact topology?
\end{Prb}

In \cite[Example 11]{Stone} it is shown that a compact space need
not be the continuous image of a compact $T_2$-space. In fact, a
careful analysis of the argument reveals the following general
fact (also stated in \cite[3.6]{Hoffmann2}).

\begin{Prp} A $KC$-space $Y$ that is the continuous image
of a compact $T_2$-space $X$ is a $T_2$-space.
\end{Prp}

{\em Proof.} Let $f:X\rightarrow Y$ be a continuous map from a
compact $T_2$-space onto a $KC$-space. Clearly $f$ is a closed
map, since $f$ is continuous, $X$ is compact and $Y$ is a
$KC$-space. The conclusion follows, since obviously  a closed
continuous image of a compact $T_2$-space, is a $T_2$-space.

\medskip
In this context also the following observation is of interest.

\begin{Prp} Let $f:X\rightarrow Y$ be a continuous map from a
maximal compact space onto a topological space $Y.$ Then $Y$ is
maximal compact if and only if the map $f$ is closed.
\end{Prp}

{\em Proof.} Suppose that $f:X\rightarrow Y$ is closed.  Since
$f^{-1}\{y\}$ is compact whenever $y\in Y$, we see that $f^{-1}K$
is compact whenever $K$ is compact in $Y$ (compare e.g. with the
proof of \cite[Theorem 3.7.2]{Engelking}). Since $f^{-1}K$ is
closed, we conclude that $K=f(f^{-1}K)$ is closed and hence $Y$ is
a compact $KC$-space. For the converse, suppose that the map
$f:X\rightarrow Y$ is not closed. Consequently there is a closed
set $F$ in $X$ such that $fF$ is not closed. Clearly the compact
set $fF$ witnesses the fact that $Y$ is not a $KC$-space.

\medskip In connection with the preceding result we note
(compare \cite[Example 3.2]{Camer}) that $T_1$-quotients of
maximal compact spaces are not necessarily maximal compact.

\begin{Prb} Are $T_1$-quotient topologies of maximal compact topologies
contained in maximal compact topologies?
\end{Prb}

Next we want to show that a weak version of our main problem has a
positive answer.

\begin{Prp} Let $(X,\tau)$ be a compact $T_1$-space. Then there is
a compact topology $\tau'$ finer than $\tau$ such that $(X,\tau')$ is a
$US$-space.
\end{Prp}

{\em Proof.} As usual two subsets $A$ and $B$ of $X$ will be
called {\em almost disjoint} provided that their intersection is
finite. Let ${\cal M}=\{A_i:i\in I\}$ be a maximal (with respect
to inclusion) family of pairwise almost disjoint injective
sequences in $X$ with a distinct $\tau$-limit (that is, each
$A_i\in {\cal M}$ is identified with $\{x_n:n\in {\bf N}\}\cup
\{x\}$ where $(x_n)_{n\in {\bf N}}$ is an injective sequence in
$(X,\tau)$ that converges to some point $x$ different from each
$x_n).$ For each $i\in I$ and $m\in {\bf N}$, let
$A_i^m=\{x_n:n\in {\bf N},n\geq m\}\cup \{x\}.$ Let $\tau'$ be the
topology on $X$ which is generated by the subbase $\tau \cup
\{X\setminus A_i^m:i\in I,m\in {\bf N}\}.$

We first show that $\tau`$ is compact. Let ${\cal C}$ be a
subcollection of ${\cal A}_\tau\cup \{A_i^n:i\in I,n\in{\bf N}\}$
with empty intersection. (Here, as in the following, ${\cal
A}_\tau$ denotes the set of $\tau$-closed sets.) Denote the
intersection of ${\cal C}$ with ${\cal A}_\tau$ by ${\cal F}.$ We
want to show that there is a finite subcollection of $\cal C$ with
an empty intersection. Of course, it will be sufficient to find a
finite subcollection of $\cal C$ with finite intersection. If
${\cal C}={\cal F},$ then such a finite subcollection of $\cal C$
must exist by compactness of $(X,\tau).$ So in this case we are
finished. If we have in our collection ${\cal C}\setminus {\cal
F}$ two sets $A_i^n$ and $A_j^m$ with $i\not=j,$ then their
intersection will be finite. So in that case we are also done.

Therefore we can assume that the set ${\cal C}\setminus {\cal F}$
is nonempty and its elements are all of the form
$A^m_{i_0}=\{x_n:n\in {\bf N},n\geq m\}\cup \{a\}$ for some fixed
$i_0\in I$ and $n\in M$ where $M$ is a nonempty subset of ${\bf N}$
and $a$ is the chosen $\tau$-limit of the sequence $(x_n)_{n\in
{\bf N}}$.

If $a\in \cap {\cal F},$ then clearly $a\in \cap {\cal C}$ ---a
contradiction to $\cap \cal {\cal C}=\emptyset.$ So there is
$F_0\in {\cal F}$ such that $a\not \in F_0.$ Since $F_0$ is
$\tau$-closed and the injective sequence $(x_n)_{n\in {\bf N}}$
$\tau$-converges to $a,$ we conclude that $F_0\cap \{x_n:n\in {\bf
N}\}$ is finite, since otherwise $a\in F_0.$ Hence for any $m\in
M$ we have that $F_0\cap A^m_{i_0}$ is finite and we are finished
again.

We deduce from Alexander's subbase theorem that the topology
$\tau'$ is compact.

Next we want to show that $(X,\tau')$ is a $US$-space. In order to
reach a contradiction, suppose that there is some sequence
$(x_n)_{n\in {\bf N}}$ that $\tau'$-converges to $x$ and $y$ where
$x$ and $y$ are distinct points in $X.$ Replacing $(x_n)_{n\in
{\bf N}}$ if necessary by a subsequence, we can and do assume that
the sequence $(x_n)_{n\in {\bf N}}$ under consideration is
injective and that $x_n$ does not belong to $\{x,y\}$ whenever
$n\in {\bf N}.$ The claim just made is an immediate consequence of
the fact that the original sequence $(x_n)_{n\in {\bf N}}$ attains
each value at most finitely many often, since $(X,\tau)$ and thus
$(X,\tau')$ is a $T_1$-space and $(x_n)_{n\in {\bf N}}$ has two
distinct limits in $(X,\tau').$

Then $(x_n)_{n\in {\bf N}}$ is an injective $\tau$-convergent
sequence having a $\tau$-limit distinct from each $x_n$ and by
maximality of the collection ${\cal M}$ there is some $A_i=\{z_n:
n\in {\bf N}\}\cup \{z\}$ where $z$ denotes the chosen
$\tau$-limit of the sequence $(z_n)_{n\in {\bf N}}$) belonging to
${\cal M}$ such that $A_i\cap \{x_n:n\in {\bf N}\}$ has infinitely
many elements. Suppose that there is some $p\in {\bf N}$ such that
$x$ or $y$ does not belong to $A_i^p$. Then $X\setminus A_i^p$ is
a $\tau'$-open neighborhood of $x$ or $y$, respectively, which
does not contain infinitely many terms of the sequence
$(x_n)_{n\in {\bf N}}$ which is impossible, because $x$ and $y$
are both $\tau'$-limits of $(x_n)_{n\in {\bf N}}.$ So there is no
such $p\in {\bf N}$ and it necessarily follows that $x=z=y$ ---a
contradiction. We conclude that $(X,\tau')$ is a $US$-space.

\begin{Cor} Each compact topology is contained in a compact $US$-topology.
\end{Cor}

\begin{Rem}
It is possible to strengthen the latter result further to the
statement that each compact topology is contained in a compact
topology with respect to which each compact countable set is
closed.
\end{Rem}

In order to see this we need the following two auxiliary results.
We recall that a topological space is called {\em sequentially
compact} provided that each of its sequences has a convergent
subsequence.

\begin{Lem}
Let $X$ be a $US$-space and let $\{K_n:n\in  {\bf N}\}$ be a
countable family of sequentially compact sets in $X$ having the
finite intersection property. Then $\cap_{n\in {\bf N}} K_n$ is
nonempty.

\end{Lem}

{\em Proof.} For each $n\in {\bf N}$ find $x_n\in \cap_{i=1}^n
K_i.$ Then the sequence $(x_n)_{n\in {\bf N}}$ has a subsequence
$(y_n)_{n\in {\bf N}}$ converging to $k\in K_1,$ because $K_1$ is
sequentially compact. Suppose that there is $m\in {\bf N}$ such
that $k\not \in K_m$. Since there is a tail of $(y_n)_{n\in {\bf
N}}$ belonging to $K_m$ and $K_m$ is sequentially compact, there
exists a subsequence of $(y_n)_{n\in {\bf N}}$ converging to some
$p\in K_m$. Since $X$ is a $US$-space, it follows that $k=p\in
K_m$
---a contradiction. We conclude that $k\in \cap_{n\in {\bf
N}}K_n.$

\begin{Lem} Each compact $US$-topology is contained in a compact
topology with respect to which each compact countable set is
closed.
\end{Lem}

{\em Proof.} Let $(X,\tau)$ be a compact $US$-space and let
$\tau'$ be the topology generated by the subbase $\tau\cup
\{X\setminus K:K\subseteq X$ is countable and compact$\}$ on $X.$

We are going to show that $\tau'$ is compact. In order to reach a
contradiction, assume that ${\cal C}$ is a subcollection of ${\cal
A}_\tau\cup \{K\subseteq X:K$ is countable and compact$\}$ having
the finite intersection property, but $\cap {\cal C}=\emptyset.$
Since $\tau$ is compact, we deduce that some compact countable set
$K$ belongs to ${\cal C}.$ Hence by countability of $K$ there must
exist a countable subcollection ${\cal D}$ of ${\cal C}$ such that
$\cap {\cal D}=\emptyset.$ Replace in ${\cal D}$ each member $F$
of ${\cal D}\cap {\cal A}_\tau$  by its trace $F\cap K$ on $K$ to
get a countable collection ${\cal D}'$ of compact countable sets
having the finite intersection property. By a result of
\cite{Levine2}, each compact countable space is sequentially
compact and hence ${\cal D}'$ is a countable collection of
sequentially compact sets in a $US$-space. Since $\cap {\cal
D}'=\emptyset,$ we have reached a contradiction to the preceding
lemma. We conclude that $\tau'$ is compact by Alexander's subbase
theorem. Evidently each compact countable set in $(X,\tau')$ is
$\tau$-compact and thus $\tau'$-closed.

\begin{Prb} Given some fixed cardinal $\kappa>\aleph_0.$
 Is each compact topology contained in a compact topology
with respect to which each compact set of cardinality $\kappa$ is
closed?
\end{Prb}

A modification of some of the arguments presented above allows us
to answer positively the variant of the main problem (see
\cite[Question 8-1, p. 56]{Cameron}) formulated for sequential
compactness instead of compactness.

\begin{Thm} Each sequentially compact topology $\tau$ on a  set $X$
is contained in a sequentially compact topology $\tau''$ that is
maximal among the sequential compact topologies on $X.$
\end{Thm}

{\em Proof.} Since $(X,\tau)$ is sequentially compact and any
convergent (sub)sequence has a constant or an injective
subsequence, it is obvious that any sequence in $(X,\tau)$ has a
subsequence that converges with respect to the supremum $\tau \vee
\tau_c$ where $\tau_c$ denotes the cofinite topology on $X.$
Therefore by replacing $\tau$ by $\tau\vee \tau_c$ if necessary,
in the following we assume that the sequentially compact topology
$\tau$ on $X$ is a $T_1$-topology.

Define now a topology $\tau'$ on $X$ in exactly the same way as
above. We next show that $(X,\tau')$ is sequentially compact
provided that $(X,\tau)$ is sequentially compact.
 Let $(y_n)_{n\in {\bf
N}}$ be any sequence in $X.$ It has a subsequence $(s_n)_{n\in
{\bf N}}$ that converges to some point $a$ in $(X,\tau)$, because
$(X,\tau)$ is sequentially compact.
 If $(s_n)_{n\in
{\bf N}}$ has a constant subsequence, then $(y_n)_{n\in {\bf N}}$
clearly has a convergent subsequence in $(X,\tau')$. So by
choosing an appropriate subsequence of $(s_n)_{n\in {\bf N}}$ if
necessary, it suffices to consider the case that $(s_n)_{n\in {\bf
N}}$ is injective and that $s_n\not=a$ whenever $n\in {\bf N}.$ By
maximality of ${\cal M}$ there is $A_i=\{z_n:n\in {\bf N}\}\cup
\{z\}$ belonging to ${\cal M}$ such that $\{s_n:n\in {\bf N}\}\cap
A_i$ is infinite. Hence there is a common injective subsequence of
the injective sequences $(s_n)_{n\in {N}}$ and $(z_n)_{n\in {\bf
N}}$ in this intersection. By definition of $\tau'$ that
subsequence converges to $z,$ because any basic
$\tau'$-neighborhood $G\cap \bigcap_{j=1}^n (X\setminus
A_j^{k_j})$ of $z$ where $G$ is $\tau$-open, $A_j\in {\cal M}$ and
$k_j\in {\bf N}$ $(j=1,\dots,n)$ contains a tail of that
subsequence, since $(z_n)_{n\in {\bf N}}$ $\tau$-converges to $z$
and $A_j\cap A_i$ is finite whenever $j=1,\dots,n.$ We conclude
that $(y_n)_{n\in {\bf N}}$ has a $\tau'$-convergent subsequence
and that $(X,\tau')$ is sequentially compact. As in the preceding
proof, one argues that $(X,\tau')$ is a $US$-space.

We now define a new topology $\tau''$ on $X$ by declaring
$A\subseteq X$ to be $\tau''$-closed if and only if $x_n\in A$
whenever $n\in {\bf N}$ and $(x_n)_{n\in{\bf N}}$ converges to $x$
in $(X,\tau')$ imply that $x\in A.$ It is well-known and readily
checked that $\tau''$ is a topology finer than $\tau'$ on $X$ with
the property that any sequence $(x_n)_{n\in {\bf N}}$ that
converges to $x$ in $(X,\tau')$ also converges to $x$ in
$(X,\tau'').$ In particular, it follows that the space
$(X,\tau'')$ is sequentially compact, because $(X,\tau')$ is
sequentially compact.

Let $K$ be a sequentially compact subset in $(X,\tau'').$ Suppose
that $x_n\in K$ whenever $n\in {\bf N}$ and that the sequence
$(x_n)_{n\in {\bf N}}$ converges to $x$ in $(X,\tau').$ Then there
is a subsequence $(y_k)_{k\in {\bf N}}$ of $(x_n)_{n\in {\bf N}}$
that converges to $r\in K$ in $(X,\tau'),$ since $K$ is
sequentially compact in $(X,\tau'')$ and $\tau'\subseteq \tau''.$
Thus $x=r$, since $(X,\tau')$ is a $US$-space and hence $x\in K.$
By the definition of the topology $\tau''$ we conclude that $K$ is
closed in $(X,\tau'').$ Therefore each sequentially compact subset
of $(X,\tau'')$ is $\tau''$-closed. By \cite[Theorem 2.4]{Camer}
we conclude that $\tau''$ is a maximal sequentially compact
topology on $X$, which is clearly finer than $\tau.$

\smallskip
Let us finally mention another possibly even more challenging
version of our main problem.

\begin{Prb} Which (compact) $T_1$-topologies are the infimum of a family of
maximal compact topologies?
\end{Prb}

Evidently the cofinite topology on an infinite set $X$ is the
infimum of the family of maximal compact Hausdorff topologies of
the one-point-compactifications $X_x$ (where $x\in X$) that we
have defined above. In Proposition 6 below we shall deal with a
special answer to Problem 5.

\section{Some further results}

Let $(X,\tau)$ be a compact topological space. Denote by ${\cal
A}_\tau$ (resp. ${\cal C}_\tau$) the set of all closed (resp.
compact) sets of $(X,\tau).$

Note that if $\tau$ and $\tau'$ are two compact topologies on a
set $X$ such that $\tau\subseteq \tau',$ then ${\cal
A}_\tau\subseteq {\cal A}_{\tau'}\subseteq {\cal
C}_{\tau'}\subseteq {\cal C}_{\tau}.$ Of course, a topology $\tau$
is a compact $KC$-topology if and only if ${\cal A}_\tau={\cal
C}_\tau.$

As usual, a collection of subsets of $X$ that is closed under
finite intersections and finite unions will be called a {\em ring
of sets} on $X.$ We consider the set ${\cal M}_\tau$ of all rings
${\cal G}$ of sets ordered by set-theoretic inclusion
 on the topological space $(X,\tau)$ such that
${\cal A}_\tau\subseteq {\cal G}\subseteq {\cal C}_\tau$. Since
${\cal A}_\tau$ is such a ring, ${\cal M}_\tau$ is nonempty. If
${\cal K}$ is a nonempty chain in ${\cal M}_\tau$, then $\bigcup
{\cal K}$ belongs to ${\cal M}_\tau$. By Zorn's lemma we conclude
that ${\cal M}_\tau$ has maximal elements.

We shall call a collection ${\cal C}$ of subsets of a set $X$ {\em
compact$^*$} provided that each subcollection of ${\cal C}$ having
the finite intersection property has nonempty intersection. We use
this nonstandard convention in order to avoid any confusion with
the concept of a compact topology.

\begin{Lem} Let $(X,\tau)$ be a compact topological space.  If
${\cal G}$ is a maximal element in ${\cal M}_\tau$ that is a
compact$^*$ collection, then ${\cal G}={\cal A}_{\tau'}$ where
$\tau'$ is a maximal compact topology finer than $\tau.$
\end{Lem}

{\em Proof.} Suppose  that ${\cal G}$ is a maximal element in
${\cal M}_{\tau}$ that is compact$^*$. Then $\{X\setminus K:K\in
{\cal G}\}$ yields the subbase of a topology ${\tau}'$ on $X.$
 Observe that ${\cal A}_\tau \subseteq {\cal
G}\subseteq {\cal A}_{{\tau}'}$. Since ${\cal G}$ is compact$^*$,
${\tau}'$ will be compact, by Alexander's subbase theorem. Because
${\tau}'$ is compact, $\tau\subseteq {\tau}'$ implies that ${\cal
A}_{{\tau}'}\subseteq {\cal C}_\tau.$ Hence ${\cal A}_{{\tau}'}$
belongs to ${\cal M}_\tau.$ We conclude that ${\cal G}={\cal
A}_{{\tau}'}$ by the maximality of ${\cal G}.$

It remains to be seen that ${\tau}'$ is maximal compact. If
${\tau}''$ is a finer topology than ${\tau}'$ and compact, then
${\cal A}_{{\tau}'}\subseteq {\cal A}_{{\tau}''}\subseteq {\cal
C}_\tau.$ Hence by maximality of ${\cal G},$ ${\cal
A}_{{\tau}''}={\cal G}={\cal A}_{\tau'}$ and so
${\tau}'={\tau}''.$ We have shown that ${\tau}'$ is maximal
compact.

\begin{Prp} Let $(X,\tau)$ be a compact topological space such that
each filterbase consisting of compact subsets has a nonempty
intersection. Then $\tau$ is contained in a maximal compact
topology ${\tau}'.$
\end{Prp}

{\em Proof.} Let ${\cal G}$ be any maximal element in ${\cal
M}_\tau$ as defined above. Recall that ${\cal G}$ is closed under
finite intersections. Hence any nonempty subcollection ${\cal G}'$
of ${\cal G}$ having the finite intersection property generates a
filterbase consisting of compact sets on $X.$ It follows from our
hypothesis that ${\cal G}$ is  a compact$^*$ collection.
Furthermore by Lemma 4 we conclude that ${\cal G}$ is equal to the
set of closed subsets of a maximal compact topology ${\tau}'$ that
is finer than $\tau.$

\medskip

It is known and easy to see (compare \cite[Theorem 6]{Levine})
that if $X$ is a compact $KC$-space, then the product $X^2$ is a
$KC$-space if and only if $X$ is a Hausdorff space. As an
application of Proposition 4 we want to show however that the
seemingly reasonable conjecture that the product topology of a
large family of maximal compact topologies is no longer contained
in a maximal compact topology is unfounded. In order to see this
we next prove the following result.

\begin{Lem} Let $(X_i)_{i\in I}$ be a nonempty family of $T_1$-spaces such
that each $X_i$ has the property that every filterbase of compact
sets has a nonempty intersection. Then the product $\Pi_{i\in I}
X_i$ also has that property.
\end{Lem}

{\em Proof.} We can (and do) assume that $I$ is equal to some
finite ordinal or an infinite limit ordinal $\epsilon.$ Let ${\cal
F}$ be a filterbase of compact subsets on the product $\Pi_{\gamma
<\epsilon}X_\gamma.$

 For each
$\alpha<\epsilon$ we shall inductively find $x_\alpha\in X_\alpha$
such that the set $A_\alpha=\{(y_\gamma)_{\gamma<\epsilon}\in
\Pi_{\gamma<\epsilon} X_\gamma:y_\gamma=x_\gamma$ whenever
$\gamma\leq \alpha\}$ satisfies $A_\alpha\cap K\not=\emptyset$
whenever  $K\in {\cal F}.$

Suppose now that for some $\delta<\epsilon$ and all
$\alpha<\delta,$ $x_\alpha\in X_\alpha$ have been chosen such that
$A_\alpha\cap K\not = \emptyset$ whenever $K\in {\cal F}.$ Let us
first establish the following claim.

Caim: $(\cap_{\alpha<\delta} A_\alpha)\cap K\not=\emptyset$
whenever $K\in {\cal F}.$

If $\delta$ is a successor ordinal, then by our induction
hypothesis $A_{\delta-1}\cap K\not =\emptyset$ whenever $K\in
{\cal F}.$ Therefore the claim is verified, since the family
$\{A_\alpha:\alpha<\delta\}$ is monotonically decreasing. So let
$\delta$ be a limit ordinal (possibly equal to $0$) and fix $K\in
{\cal F}.$ Since for each $\alpha<\delta,$ $A_\alpha$ is closed
because every space $X_\alpha$ is a $T_1$-space, and since
$A_\alpha\cap K\not=\emptyset$  the claim holds by compactness of
$K$ and the monotonicity of the sequence
$\{A_\alpha:\alpha<\delta\}$. (For the case that $\delta=0$ as
usual we use  the convention that
$\cap\emptyset=\Pi_{\gamma<\epsilon}X_\gamma.$)

\medskip
Continuing now with the proof we next consider the filterbase
$\{[\mbox{pr}_{X_\delta}(\cap_{\alpha<\delta}A_\alpha\cap K)]:K\in
{\cal F}\}$ of compact sets on $X_\delta.$ By our assumption on
$X_\delta$ there exists some $$x_\delta\in \cap_{K\in {\cal F}}
[\mbox{pr}_{X_\delta}(\cap_{\alpha<\delta}A_\alpha\cap K)].$$ It
remains to show that for each $K\in {\cal F},$
$\{(y_\gamma)_{\gamma<\epsilon}\in
\Pi_{\gamma<\epsilon}X_\gamma:y_\gamma=x_\gamma,\gamma\leq
\delta\}\cap K\not=\emptyset;$ but this is an immediate
consequence of $x_\delta \in
\mbox{pr}_{X_\delta}(\cap_{\alpha<\delta}A_\alpha\cap K).$ Finally
note that $\cap_{\alpha
<\epsilon}A_\alpha=\{(x_\alpha)_{\alpha<\epsilon}\}$ and that ---
for $\epsilon$ exactly as in the case of the ordinal $\delta$
above
--- $\cap_{\alpha<\epsilon}A_\alpha\cap K\not=\emptyset$ whenever
$K\in {\cal F}.$ Hence the assertion of the lemma holds.

\medskip
\begin{Prp} The product topology of a nonempty family of compact $KC$-topologies is
 contained in a maximal compact topology.
\end{Prp}

{\em Proof.} Note first that in a compact $KC$-topology each
filterbase of compact sets has a nonempty intersection.
 We conclude by the preceding lemma and Proposition 4 that the compact product
topology of an arbitrary nonempty family of maximal compact
topologies is contained in a maximal compact topology.

\begin{Cor} Let $(X_i)_{i\in I}$ be a nonempty family of spaces
each of which is contained in a maximal compact topology. Then
also their product topology is contained in a maximal compact
topology.
\end{Cor}

\section{Sobriety and maximal compactness}

\smallskip
Note that the closure of each irreducible subspace of a
topological space is irreducible. Recall also that a topological
space is called {\em sober} (see e.g. \cite{Hoffmann}) provided
that every irreducible closed set is the closure of some unique
singleton. Clearly each Hausdorff space is sober. Furthermore a
subset of a topological space is called {\em saturated} provided
that it is equal to the intersection of its open supersets.

A short proof of the following result is given in
\cite{KeimelPaseka}.

Let $\{K_i: i\in I\}$ be a filterbase of (nonempty) compact
saturated subsets of a sober space $X.$ Then $\bigcap_{i\in I}K_i$
is nonempty, compact, and saturated, too; and an open set $U$
contains $\bigcap_{i\in I} K_i$ iff $U$ contains $K_i$ for some
$i\in I.$

\begin{Cor} Let $(X,\tau)$ be a compact sober $T_1$-space.
Then $\tau$ is contained in some maximal compact topology
${\tau}'.$
\end{Cor}

{\em Proof.} Since all (compact) sets in a $T_1$-space are
saturated, the condition stated in Proposition 4 is satisfied by
the result just cited. The statement then follows from Proposition
4.

\begin{Prb} Characterize those sober compact topologies that are
contained in a maximal compact topology.
\end{Prb}

\begin{Rem} {\em Let us observe that the maximal compact
topology $\tau'$ obtained in Corollary 3 will be sober, since the
only irreducible sets with respect to the coarser topology $\tau$
are the singletons. Van Douwen's example \cite{vanDouwen}
mentioned earlier shows that a maximal compact topology need not
be (contained in) a compact sober topology. }
\end{Rem}

\begin{Ex} {\em Note that the closed irreducible subsets of the
one-point-compactification $X$ (of the Hausdorff space) of the
rationals are the singletons: Any finite subset of a $T_1$-space
with at least two points is discrete and hence not irreducible.
Moreover any infinite subset of $X$ contains two distinct
rationals and thus cannot be irreducible.
 We conclude that an
arbitrary power of $X$ is a compact, sober $T_1$-space, because
products of sober spaces are sober (see e.g. \cite[Theorem 1.4]
{Hoffmann}). }
\end{Ex}

In the light of the proof of Proposition 4 one wonders which
compact sober $T_1$-topologies can be represented as the infimum
of a family of maximal compact topologies (compare Problem 5). Our
next result provides a partial answer to this question. An
interesting space satisfying the hypothesis of Proposition 6 is a
$T_1$-space constructed in \cite{KunziWatson}: It has infinitely
many isolated points although each open set is the intersection of
two compact open sets. (It was noted in the discussion \cite[p.
212]{KunziWatson} that that space is locally compact and sober.)

\smallskip

Recall that a topological space $X$ is called {\em sequential}
 (see \cite[p. 53]{Engelking}) provided that
a set $A\subseteq X$ is closed if and only if together with any
sequence it contains all its limits in $X$.

\begin{Prp} Each compact sober  $T_1$-space $(X,\tau)$ which is locally compact or
sequential
 is the infimum of a family of maximal compact
topologies.
\end{Prp}

{\em Proof.} Note that if $K$ belongs to the closed sets of a
maximal compact topology $\sigma$ finer than $\tau,$ then $K$ is
compact with respect to $\sigma$ and thus with respect to $\tau$.
In order to verify the statement, it therefore suffices to
construct for any compact set $C$ that is not closed in $(X,\tau)$
a maximal compact topology $\sigma$ finer than $\tau$ in which $C$
is not closed.

So let $C$ be a compact set that is not closed in $(X,\tau).$ In
$(X,\tau)$ we shall next find a compact set $K_0$ such that
$K_0\cap C$ is not compact.

Suppose first that $X$ is locally compact.

Then
there is $x\in X$ such that $x\in \mbox{cl}_\tau{C}\setminus C.$
Let ${\cal F}=\{K:K$ is a compact neighborhood at $x$ in
$(X,\tau)\}.$ Of course, $\cap {\cal F}=\{x\},$ since $X$ is a
locally compact $T_1$-space. Suppose that $K\cap C$ is compact in
$(X,\tau)$ whenever $K\in {\cal F}.$

Then $\{K\cap C: K\in {\cal F}\}$ is a filterbase of compact
saturated sets in $X.$ According to the result cited above from
\cite{KeimelPaseka}, we have $\cap {\cal F}\cap C=\bigcap \{K\cap
C:K\in {\cal F}\}\not=\emptyset.$ Since $x\not\in C,$ we have
reached a contradiction. Thus there is a compact neighborhood
$K_0$ of $x$ such that $K_0\cap C$ is not compact in $(X,\tau)$.

Suppose next that $X$ is sequential. Since $C$ is not closed,
there is a sequence $(x_n)_{n\in {\bf N}}$ in $C$ converging to
some point $x\in X$ such that $x$ does not belong to $C.$ Assume
that $\{(\{x\}\cup \{x_n:n\in {\bf N},n\geq m\})\cap C:m\in {\bf
N}\}$ is a filterbase of compact sets. Clearly its intersection is
empty, because $\tau$ is a $T_1$-topology and $(x_n)_{n\in {\bf
N}}$ converges to $x$---a contradiction. Hence there is $m\in {\bf
N}$ such that $(\{x\}\cup \{x_n:n\in {\bf N},n\geq m\} )\cap C$ is
not compact. Denote the compact set $\{x\}\cup \{x_n:n\in {\bf
N},n\geq m\}$
  by $K_0.$

   So our claim holds in either case.

Note now that $\tau\cup \{X\setminus K_0\}$ is a subbase for a
compact topology $\tau'$ on $X$ that is also sober and $T_1.$ By
Corollary 3 there is a maximal compact topology $\tau''$ finer
than $\tau'.$ Observe that $X\setminus C\not\in \tau'':$ Otherwise
$C\in {\cal A}_{\tau''}$ and, since $K_0\in {\cal A}_{\tau''},$
also $K_0\cap C \in {\cal A}_{\tau''}.$ Therefore $K_0\cap C\in
{\cal C}_{\tau''}$ and $K_0\cap C\in {\cal C}_{\tau}$ ---a
contradiction. Thus indeed $X\setminus C\not\in \tau''.$ We
conclude that $\tau$ is the infimum of a family of maximal compact
topologies.

\medskip
Observe that the argument above also yields the following results.

\begin{Cor} Each locally compact  sober $T_1$-space in which
the intersection of any two compact sets is compact is a $KC$-space
(and therefore is a regular Hausdorff space).
\end{Cor}

\begin{Cor} Each sequential sober $T_1$-space in which the intersection
of any two compact sets is compact is a $KC$-space.
\end{Cor}

\smallskip
We next give an example of a compact sober $T_1$-topology that is
not the infimum of a family of maximal compact topologies.

\begin{Ex}
{\em Let $Y$ be an uncountable set and let $-\infty$ and $\infty$
be two distinct points not in $Y$. Set $X=Y\cup
\{-\infty,\infty\}.$ Each point of $Y$ is supposed to be isolated.
The neighborhoods of $\infty$ are the cofinite sets containing
$\infty$ and the neighborhoods of $-\infty$ are the cocountable
sets containing $-\infty.$ Clearly $X$ is a compact sober
$T_1$-space.

 Next we show that with respect to the defined topology
$\tau$ a subset $A$ of $X$ is compact and not closed if and only
if $A$ is uncountable, $\infty\in A$ and $-\infty\not\in A:$
Indeed, if $\infty\in A,$ then $A$ is clearly compact and if $A$
is uncountable and $-\infty\not \in A,$ then $A$ cannot be closed.
In order to prove the converse suppose that $A$ is compact and not
closed in $(X,\tau).$ Then $A$ is certainly infinite. It therefore
follows from compactness of $A$ that $\infty\in A.$ Since $A$ is
not closed, we conclude that $-\infty\in \overline{A}\setminus A$
and hence $A$ is uncountable.

Of course, if $\tau'$ is a maximal compact topology such that
$\tau \subseteq \tau'$, then ${\cal A}_\tau \subseteq {\cal
A}_{\tau'}\subseteq {\cal C}_{\tau'}\subseteq {\cal C}_\tau.$
Observe that the topology $\tau''$ generated by the subbase
$\{\{-\infty\}\}\cup \tau$ clearly yields a compact $T_2$-topology
finer than $\tau.$ Obviously, ${\cal C}_\tau\setminus {\cal
A}_\tau\subseteq {\cal A}_{\tau''}$ by the description found above
of the nonclosed compact sets in $(X,\tau).$ Thus ${\cal
A}_{\tau''}={\cal C}_\tau.$ We conclude that $\tau''$ is finer
than any maximal compact topology containing $\tau.$ Hence
$\tau''$ is the only maximal compact topology (strictly) finer
than $\tau.$}

\end{Ex}

\medskip

Let us recall that a topological space is called {\em strongly
sober} provided that the set of limits of each ultrafilter is
equal to the closure of some unique singleton. Of course, each
compact Hausdorff space satisfies this condition.

We finally observe that  each locally compact strongly sober
topological space $(X,\tau)$ possesses a finer compact Hausdorff
topology; just take the supremum of $\tau$ and its dual topology
(see e.g. \cite[Theorem 4.11]{Kopperman}). By definition, the
latter topology is generated by the subbase $\{X\setminus K: K$ is
compact and saturated in $X \}$ on $X.$

No characterization seems to be known of those topologies that
possess a finer compact Hausdorff topology.

\bigskip
 \noindent
Dept. Math. Appl. Math. \\
University of Cape Town \\
Rondebosch 7701 \\
South Africa \\
kunzi@maths.uct.ac.za

\bigskip
\noindent Math. Institute\\
University of Berne \\
Sidlerstrasse 5\\
3012 Berne\\
Switzerland\\
dominic.vanderzypen@math-stat.unibe.ch

 \end{document}